\documentclass{amsart}

\author{Gary Froyland}
\address{School of Mathematics and Statistics,
University of New South Wales, Australia}
\email{G.Froyland@unsw.edu.au}

\author{Georg A. Gottwald}
\address{School of Mathematics and Statistics,
University of Sydney, Australia}
\email{georg.gottwald@sydney.edu.au}

\author{Andy Hammerlindl}
\address{School of Mathematics and Statistics,
University of Sydney, Australia
and School of Mathematics and Statistics,
University of New South Wales, Australia}
\email{andy@maths.usyd.edu.au}

\title{A trajectory-free framework for analysing multiscale systems}
\thanks{This project was funded by the Australian Research Council Grant DP120104514.}

\usepackage{amsmath}
\usepackage{amssymb}
\usepackage{amsfonts}
\usepackage{amsthm}
\usepackage{bbm}
\usepackage{bm}
\usepackage{graphicx}
\usepackage{setspace}
\usepackage{booktabs}

\newcounter{algonum}

\newcommand{\bbR}{\mathbb{R}}

\newcommand{\inv}{^{-1}}

\newcommand{\cX}{\mathcal{X}}
\newcommand{\cZ}{\mathcal{Z}}

\newcommand{\cL}{\mathcal{L}}

\newcommand{\qandq}{\quad\text{and}\quad}

\newcommand{\avg}{_{\operatorname{avg} }}
\newcommand{\err}{\operatorname{err} }

\newcommand{\ep}{\epsilon}
\newcommand{\lam}{\lambda}

\newcommand{\sig}{\sigma}

\newcommand{\cP}{\mathcal{P}}

\newcommand{\bbE}{\mathbb{E}}

\newcommand{\Real}{\operatorname{Re}}

\newcommand{\ideallam}{{-k}/{\ep}}

\newcommand{\muA}{\mu^A}
\newcommand{\DA}{D^A}

\newcommand{\mun}{\mu^{\text{nor}}}
\newcommand{\mut}{\mu^{\text{tan}}}
\newcommand{\DDn}{D^{\text{nor}}}
\newcommand{\DDt}{D^{\text{tan}}}

\newcommand{\cLa}{\cL^*}

\newcommand{\cLtan}{{\widehat \cL}_{\text{tan}}}
\newcommand{\lamtan}{\widehat \lam}

\newcommand{\Ito}{It\=o}

\theoremstyle{remark}

\newcommand{\Nmat}{\ensuremath{50}}
\newcommand{\Lmat}{\ensuremath{5}}
\newcommand{\lamone}{\ensuremath{2.0971 \times 10^{-11}}}
\newcommand{\lammat}{\ensuremath{-0.6467 + 0.1097 i}}
\newcommand{\avgdataone}{\ensuremath{1.4809 \times 10^{3}}}

\newcommand{\avgdatathree}{\ensuremath{1.4216}}
\newcommand{\avgdatafour}{\ensuremath{704.27}}
\newcommand{\avgdatafive}{\ensuremath{25.165}}

\newcommand{\minangle}{\ensuremath{55}}
\newcommand{\maxangle}{\ensuremath{125}}

\begin{document} 

\begin{abstract}
    We develop algorithms built around properties of the 
    transfer operator and Koopman operator which
    1) test for possible multiscale dynamics in a given
    dynamical system,
    2) estimate the magnitude of the time-scale separation,
    and finally
    3) distill the reduced slow dynamics on a suitably designed
    subspace. 
    By avoiding trajectory integration,
    the developed techniques
    are highly computationally efficient.
    We corroborate our findings with numerical
    simulations of a test problem. 
\end{abstract}
\maketitle

\section{Introduction} \label{sec-intro} 

Effective numerical simulation of multiscale systems constitutes a formidable challenge. 
Consider a system which has slow dynamics on a time-scale of order one and fast
dynamics on the scale of order $1/\ep$ for some
parameter $\ep \ll 1$.
In molecular dynamics, for example, the time-scale separation can be of the order of $\ep\approx 10^{-14}$ \cite{E11}.
To accurately simulate orbits numerically and to assure numerical stability,
the time step of the integrator must be of the order of $\ep$.
To capture the relevant slow dynamics a total number of integration steps of
the order of $1/\ep$ is required, making direct numerical simulations of
orbits computationally impractical.

Numerical integrators are subject to two main sources of error.
The first is truncation error, which is
the inability of the numerical method
(Runge--Kutta, Euler--Maruyama, etc.) to fully capture the actual dynamics
of the system.
The second is round-off error, due to implementing
the numerical method on a computer with finite precision arithmetic.
While truncation error decreases with a smaller time step,
round-off error increases
\cite{HenriciBook, HairerEtAlBookI}.
In a multiscale system, if the time-scale separation is large, it may be
impossible to find a time step which is simultaneously small enough to avoid
significant truncation error for the fast dynamics and 
sufficiently large to avoid detrimental accumulation
of round-off error for the slow
dynamics.


Even if orbits could be computed exactly, analysing a multiscale system using a
time series extracted from a true orbit can still yield incorrect data about
the diffusion process of the slow variables
\cite{PavliotisStuart2007}.
To avoid this problem, the time series must be sampled at a rate intermediate
between the slow and fast variables
and these rates might not be known in advance.

There exists a variety of numerical methods dealing with one or more aspects of these numerical difficulties (see \cite{E11} and references therein). These methods rely on producing trajectories of the dynamical system via some form of time-integration with some of the issues mentioned above remaining. In this paper, we develop algorithms which avoid trajectory integration altogether. Besides the advantages relating to the issues of time-integration mentioned above, the algorithm allows for a huge reduction in computational time. Our main objective is to develop numerical algorithms which, given a dynamical system,
\begin{enumerate}
\item test whether the system exhibits multiscale behaviour, and if so
\item determine the order of the time-scale separation, and then 
\item construct effective reduced equations for the slow dynamics allowing for the application of large time steps. 
\end{enumerate}
The framework we adopt for this integration-free approach is based on the infinitesimal generator associated to the underlying continuous-time dynamical system. In Section~\ref{sec-generator} we briefly review the notion of generators of transfer and Koopman operators. Section~\ref{sec-algo} introduces a trajectory-free test for multiscale behaviour. The degree of time-scale separation is estimated in algorithms described in Sections~\ref{sec-timescale} and \ref{sec-arc}. A method to determine the reduced slow dynamics from a multiscale system without relying on statistics obtained from long time-integrations is given in Section~\ref{sec-reduced}. The algorithms are tested in numerical simulations in Section~\ref{sec-numerics}. We conclude with a discussion in Section~\ref{sec-discussion}.

\section{Generators} \label{sec-generator} 
We describe our methodology for \Ito\ drift-diffusion processes, as these are a large and flexible class of dynamical systems, and the spectral properties of the corresponding transfer operators are relatively straightforward.
Consider a drift-diffusion process
\begin{equation} \label{eqn-dd-process}
    d \zeta_i = \mu_i\,dt + \sum_{k=1}^\ell \sig_{ik}\, dW_k \quad
    \text{with $i=1,\ldots,d$}
\end{equation}
defined on a subset $\cZ$
of $\bbR^d$ where $l  \le  d$ and
each $W_k$ for $k=1,\ldots,\ell$
represents an independent Wiener process.
Given a probability density function at time $t=0$,
the density at future times
is determined by the
Fokker--Planck
equation
\[
    \frac{\partial \rho}{\partial t} = \cL \rho
\]
where
\begin{equation} \label{eqn-FP}
    \cL \rho = - \sum_{i=1}^N
             \frac{\partial}{\partial z_i} \bigl[\mu_i \,\rho \bigr]
             + \frac{1}{2}\sum_{i,j=1}^N
             \frac{\partial^2}{\partial z_i\,\partial z_j}
             \bigl[D_{ij}\, \rho \bigr].
\end{equation}
The second order differential operator $\cL$ is called the
{\em Fokker--Planck operator}.
At each point $z$ in the phase space $\cZ$,
the vector $\mu(z) \in \bbR^d$ represents the drift of the process
and the positive semi-definite matrix
$D(z) = \sig(z) \sig(z)^\top \in \bbR^{d \times d}$
the diffusion.
The operator $\cL$ generates a family of operators
$e^{t \cL}$ for $t \ge 0$ such that
$e^{s \cL}e^{t \cL}=e^{(s+t) \cL}$
and
$\lim_{t \to 0} \tfrac{1}{t}[e^{t \cL} - Id] = \cL$.
If $\rho \in L^1(\cZ)$ is an initial probability density, then
$e^{t \cL}(\rho)$
is the density after time $t$.
Thus, $e^{t \cL}$ may be thought of as a
transfer operator
defined on $L^1(\cZ)$.
Since the underlying system is a random dynamical system
and $e^{t \cL}$ represents an average over all possible random paths,
it is an {annealed} transfer operator
\cite{Baladi1997}.

In many cases, the Fokker--Planck operator has compact
resolvent \cite{Helffer2005}.
In particular, its spectrum consists of a countable set of eigenvalues
$\{\lam_k\}_{k=0}^\infty$ which, when ordered by the convention
\begin{equation} \label{eqn-decreal}
    0 = \Real \lam_0  \ge  \Real \lam_1  \ge  \Real \lam_2  \ge  \cdots,
\end{equation}
satisfy $\lim_{k \to \infty} \Real \lam_k \to -\infty$.
As a consequence,
for each $t>0$ the operator $e^{t \cL}$ is compact
with eigenvalues $e^{t \lam_k}$ tending to zero as $k \to \infty$.
The invariant density $\rho_0$ of the system
is an eigenfunction of $\cL$ associated to the eigenvalue $\lam_0 = 0$.

%
%

\medskip

The Kolmogorov backward equation is given by
$\tfrac{\partial f} {\partial t} = \cLa f$
where the adjoint
of the Fokker--Planck operator is given by
\begin{equation} \label{eqn-FPa}
    \cLa f = \sum_{i=1}^N
             \mu_i \frac{\partial f}{\partial z_i}
             + \frac{1}{2}\sum_{i,j=1}^N
             D_{ij} \frac{\partial^2 f}{\partial z_i\,\partial z_j}.
\end{equation}
This adjoint operator generates a family of operators
$K_t := e^{t \cLa} = (e^{t \cL})^*$
for $t  \ge  0$.
If $f \in L^\infty(\cZ)$,
then $K_t(f) \in L^\infty(\cZ)$
is given by
$K_t(f)(z) = \bbE f(\zeta(t))$
where the expectation is over all paths $\zeta(t)$ in 
the drift-diffusion process which satisfy $\zeta(0)=z$.
The operator
$K_t$ may therefore be regarded as
an annealed {Koopman} or composition operator.


The operators $\cLa$ and $\cL$ share the same eigenvalues $\lam_k$. Consider the eigenfunction $\psi_k \in L^\infty(\cZ)$ of 
$\cLa$ associated to $\lam_k$.
This function $\psi_k$
is an observable which evolves according to
\[
    K_t \psi_k = e^{t \lam_k} \psi_k
\]
and therefore decays to zero at the rate given by
$|e^{t \lam_k}|$ for $t > 0$.
In general, the eigenfunctions of $\cLa$
associated to eigenvalues with real part closest to zero
are the observables
of the system which decay at the slowest speeds possible.
The eigenvalues $\lam_k$ for small $k$ therefore
correspond to speed of the slowest dynamics present in the system.
These ideas have been exploited to identify almost-invariant sets
\cite{dellnitz1999,fjk2013,FGTQ2014} for deterministic and annealed random dynamics, and coherent sets \cite{FLS2010} for quenched random or time-dependent dynamics.

Suppose that the system under study is a multiscale system.
That is, assume there is a value $0 < \ep \ll 1$ representing the time-scale
separation
and transverse directions of fast and slow dynamics such that
the components of the drift and diffusion
are on the order of $O(\ep \inv)$ in the fast direction
and $O(1)$ in the slow direction.
Further suppose
that $\lam_k$ is an eigenvalue of $\cLa$ such that $|\Real \lam_k| \ll \ep \inv$.
Then, one may choose an intermediate time $t \gg \ep$
such that
\[
    |\Real t \lam_k| \ll 1   \quad \Rightarrow \quad   e^{t \lam_k} \approx 1   \quad \Rightarrow \quad   K_t \psi_k \approx \psi_k.
\]
Since $t \gg \ep$, the system after time $t$ has evolved towards
equilibrium in the fast dynamics implying that $K_t \psi_k$
is nearly constant along the direction of the fast dynamics.
Therefore, $\psi_k$ should be approximately constant in this direction
as well.
This observation is developed in more detail in
\cite{Crommelin11,FGH1}.

Given the ordering in \eqref{eqn-decreal},
we refer to the eigenvalues $\lam_k$ for small $k$
(and therefore small $|\Real \lam_k|$)
as the \emph{leading} eigenvalues
with the associated \emph{leading} eigenfunctions $\psi_k$ .
By the above observations, in a multiscale system,
the level sets of a leading eigenfunction $\psi_k$ closely approximate
\emph{fast fibres}, submanifolds of the system along which the fast dynamics
occurs.
This property is used in the following sections to develop algorithms
which test for multiscale behaviour.


\section{An algorithm to test for multiscale behaviour} \label{sec-algo} 



%
%
%
%
%
%
%

We now use these properties of the Fokker--Planck operator and its adjoint to
develop an algorithm to test for multiscale behaviour.
We focus on the case of 
stochastic differential equations (SDEs)
where the slow dynamics is
one-dimensional.
The only information used by the algorithm are functions
$\mu$ and $D$ defining the
drift and diffusion of the process.  In particular, we do not assume
any a priori knowledge of the slow or fast directions (if they exist) or that
these directions align with the coordinates in which the system is defined.
We first outline the idea behind the algorithm before providing the actual algorithm.


As established in the previous section, the level sets of the leading eigenfunctions $\psi_k$ of the adjoint Fokker-Planck operator $\cLa$ contain information about the multi-scale behaviour of the dynamics. If the system has multiscale behaviour,
these level sets will approximate a fast fibre of the system. We use the eigenfunction $\psi_1$ and numerically compute a fibre $F$ as a connected component of a level set $\psi_1 \inv(\{s\})$ (cf.~step 3 in Algorithm~1). If the leading eigenfunctions are complex valued, we may choose without loss of generality the real (or imaginary) part of the eigenfunctions to define the fibres (in Section \ref{sec-numerics} we will give such an example). There is a significant history of computing eigenfunctions of the Fokker--Planck operator or its adjoint numerically \cite{Boyd01, fjk2013} which we will employ to determine  eigenfunctions $\psi_k$ of $\cLa$  associated to the eigenvalues $\lam_k$ for small $k$ (cf. step 2 in Algorithm~1).

If the system actually exhibits multiscale dynamics and the fibre $F$ approximates a fast fibre, then the speed of the dynamics will be much faster along the fibre than transverse to it. Therefore, one expects in this case that the drift and the diffusion along one or both of the averages tangent to the fibre $F$ will be significantly larger than those corresponding to the direction normal to $F$. We therefore calculate the drift and diffusion coefficients for a set of points on the fibre $F$ and then subsequently use their averages along the fibre to probe for fast dynamics along the fibre and slow dynamics transverse to the fibre (this is achieved in steps 4 and 5 of Algorithm~1, respectively).

Computing the drift and diffusion coefficients $\mun(z)$, $\mut(z)$, $\DDn(z)$, and $\DDt(z)$ normal and tangent to the fibre at a finite set of points $Q$ uniformly distributed over $F$ constitutes the most complicated step of the algorithm and is explained in further detail below in Subsection~\ref{sec-normtan}.

The most natural way to compute the averages of the drift and diffusion coefficients $\mun(z)$, $\mut(z)$, $\DDn(z)$, and $\DDt(z)$ along the fibre is using the invariant density along the fibre. As an approximation for the invariant density of the fast dynamics on the fibre, we use the invariant density function $\rho_0$ of the full multi-scale system and simply restrict $\rho_0$ to $F$. The property that $\rho_0$ is invariant is equivalent to $\cL \rho_0 = 0$. Therefore, we estimate $\rho_0$ by solving numerically for the eigenfunction of $\cL$ whose associated eigenvalue is closest to zero (cf. step 1 of Algorithm~1). The measure on $F$ can be approximated by a measure supported on a finite set. We choose a finite set of points $Q$ uniformly distributed over $F$ and define weights $w_z \in \bbR$ such that $\sum_{z \in Q} w_z \delta_z$ approximates the density $\rho_0$ on $F$. Here, $\delta_z$ is the Dirac measure at $z$ and $w_z := c \rho_0(z)$ where the scaling constant $c$ is defined such that $\sum_{z \in Q} w_z = 1$ (cf. step 4 of Algorithm~1).

As a criterion for the presence or absence of multi-scale behaviour we then compare the averages along the fibre of the absolute values of the normal and tangential drift and diffusion coefficients; if the tangential components are much larger than the perspective normal ones, this is indicative of $F$ being a fast fibre and the system exhibiting multi-scale behaviour (cf. step 7).

\subsection{Computing normal and tangent components}
\label{sec-normtan}
We now explain step 5 of Algorithm 1 in detail.
To compute the values
$\mun(z)$, $\mut(z)$, $\DDn(z)$, and $\DDt(z)$
one must split the dynamics at a point $z \in F$
into the dynamics along the fibre $F$
and dynamics normal to $F$.
There is a subtle issue caused by the diffusion.
If the embedding of the fibre into the overall phase space is non-linear,
then diffusion along the fibre will generate a drift term for the full system.
As a simple example of this,
consider the embedding of $\bbR$ into $\bbR^2$ given by
$u \mapsto (u,u^2)$
and consider a (drift-free) diffusion process $du = dW$
where $W$ is the standard Wiener process on $\bbR$.
For a point $(x,y)$ on this embedded curve,
\Ito's lemma \cite{OksendalBook} gives
\begin{align*}
    dx &= du = dW, \\
    dy &= d (u^2)
       = 2 u\, du + \frac{1}{2} 2\, du^2
       = 2 x\, dW + dt.
\end{align*}
In the notation of \eqref{eqn-dd-process},
this corresponds to a drift vector $\mu = (0,1)$ at every point
and a diffusion matrix $D = \sig \sig^\top$ with
$\sig = (1, 2x)^\top$
at each point $(x,y)$ on the parabola.
Notice that while $\sig$ is everywhere tangent to the parabola,
the drift vector $\mu$ always points in the vertical $y$ direction
and is not tangent to the parabola.
To give an intuition for this, one can consider a point mass at
$(x,y) = (0,0)$ at time $t=0$.
As time progresses, the mass spreads out along the parabola
$\{(x,y) : y=x^2\}$ and at each time $t > 0$
the center of mass (as calculated in $\bbR^2$)
is at the point $(0,t)$ off of the parabola.
Because of this effect of diffusion along a fibre
leading to drift for the full system,
we must apply a non-linear change of coordinates to determine
the components of the drift and diffusion both normal and tangent
to the fibre.

\begin{figure}
    [tp]
    \begin{centering}
        \includegraphics{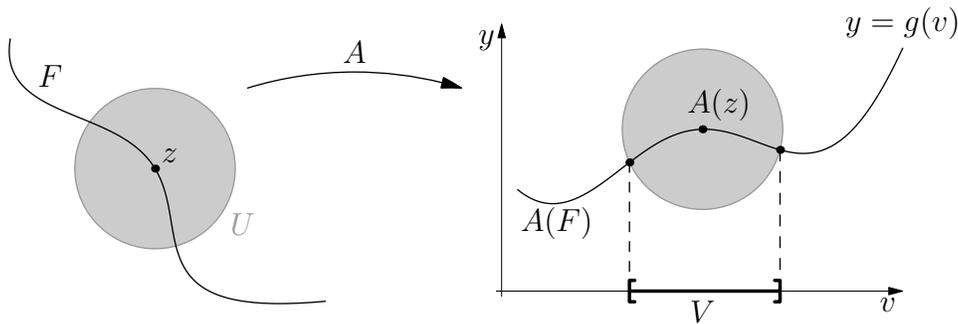}  \end{centering}
    \caption{
    Once a fast fibre $F$ is identified,
    a linear transformation $A$ is applied so that $A(F)$
    in a neighbourhood of $A(z)$ may be expressed
    as the graph $\{(v,g(v)) : v \in V\}$ of a function
    $g: V \to \bbR$ where $V \subset \bbR^{d-1}$.
    A non-linear transformation can then be applied which takes
    $A(F)$ to the flat subset $V \times \{0\} \subset \bbR^d$.
    In these coordinates, the tangent and normal
    components of the dynamics can be analysed directly.
    }
    \label{fig-vygraph}  \end{figure}

%
This may be done in two steps as depicted in
Figure \ref{fig-vygraph}.
Let $U \subset \cZ$
be a neighbourhood of the point $z \in F$.
First, find a linear change of coordinates $A:\bbR^d \to \bbR^d$,
given by a $d \times d$ matrix,
so that $A(F \cap U)$ is equal to the graph
$
    \{ (v, g(v)) : v \in V \}
$
of a function $g:V \to \bbR$ for some $V \subset \bbR^{d-1}$.
As $F$ is the level set of an eigenfunction,
such a graph exists due to the implicit function
theorem
so long as the gradient of the eigenfunction is non-zero at $z$.

A subsequent non-linear change of coordinates $(v,y) \mapsto (v,y-g(v))$
flattens the graph out,
mapping $A(F)$ to the subset $V \times \{0\} \subset \bbR^{d-1} \times \bbR$.
After transforming the dynamics to lie in the hyperplane
$\bbR^{d-1} \times \{0\}$,
it is straightforward to isolate the normal and tangent components of the
diffusion process.

For completeness,
we give here the formulas for transforming the drift and diffusion.
These formulas can be derived either by \Ito\ calculus on
\eqref{eqn-dd-process} or by applying a deterministic change of variables
to the linear operator $\cLa$. 
For a linear change of coordinates $A:\bbR^d \to \bbR^d$
the drift vector $\mu(z)$ at $z$ is mapped to $A \mu(z)$ now at the point $A(z)$
and
the diffusion matrix $D(z)$ is mapped to $A D(z) A^\top$.
In what follows, we write $\muA$ and $\DA$ for
drift and diffusion functions after this linear transformation.

Now consider the change of coordinates $(v,y) \mapsto (v,y-g(v))$
which in all $d$ coordinates
may be written as
\[
    (z_1, \ldots, z_{d-1},\,  z_d) \mapsto
    (z_1, \ldots, z_{d-1},\, z_d - g(z_1,\ldots,z_{d-1})).
\]
This transforms the drift-diffusion process given by $\muA$ and $\DA$
to one given by $\hat \mu$ and $\hat D$
where
\begin{align} \label{eqn-hat-mu}
    \hat \mu_i &= \muA_i \quad &&\text{for $i<d$,}\\
    \hat \mu_d &= \muA_d - \sum_{i<d} \muA_i g_i +
    \frac{1}{2} \sum_{i,j<d} \DA_{ij} g_{ij}
\end{align}
and

\begin{align} \label{eqn-hat-D}
    \hat D_{ij} &= \DA_{ij} \quad &&\text{for $i,j<d$},\\
    \hat D_{id} &= \DA_{id} - \sum_{j<d} \DA_{ij} g_j \quad &&\text{for $i<d$, and}\\
    \hat D_{dd} &=
    \DA_{dd} - 2 \sum_{i<d} \DA_{id} g_i + \sum_{i,j<d} \DA_{ij} g_i
    g_j.\label{eqn-last-D}
\end{align}
In the above equations, $g_i$ and $g_{ij}$ are the partial derivatives of $g$ with respect to $x$.
Since only first and second-order partial derivatives appear above,
only a second-order approximation of $g$ is necessary
in order to compute the new values. 

We propose the following method to realize these
two changes of variables numerically.
Suppose the fibre $F$ in a neighbourhood of a point $z$ is represented
numerically by a collection of points $\{q_n\}$ near $z$.
First, perform linear regression on this set of points
to get a hyperplane  $H \subset \bbR^d$
such that the set $H+z$ is tangent to $F$ at $z$.
Using a numerical orthogonalization process such as the QR decomposition,
construct an orthonormal basis $\{u_1, \ldots, u_{d-1}, u_d\}$ of $\bbR^d$
such that $u_1, \ldots, u_{d-1} \in H$.
This basis then gives the rows of an orthogonal $d \times d$ matrix
representing $A$.

Next, write each point $A(q_n)$
in the form $(v_n, y_n)$ with
$v_n \in \bbR^{d-1}$
and $y_n \in \bbR$.\footnote{Here and later in the paper,
we adopt the convention that subscripts with the letters $m$ and $n$
refer to distinct points in the phase space $\cZ\subset\bbR^d$ and
subscripts with $i$ and $j$ denote the coordinates of a point.
That is, $z_i,z_j \in \bbR$ for $1 \le i,j \le d$
would refer to the coordinates of a point $z = (z_1,\ldots,z_d)\in \cZ$
whereas
$q_m$ and $q_n$
would refer to points in $\cZ$.
Throughout the paper, the subscript $k$ is used only to index
the eigenvalues and eigenfunctions of a linear operator.
}
By polynomial regression or a similar fitting technique,
find a polynomial $p_g:\bbR^{d-1} \to \bbR$
such that $p_g(v_n) \approx y_n$ and where this approximation
is to at least second order.
Using $p_g$ in place of $g$, compute $\hat \mu$
and $\hat D$ as above.


If the fibre $F$ truly corresponds to a fast fibre in a multiscale system,
one expects the drift and diffusion terms to be large for coordinates
tangent to $F$ and small otherwise.
As a test of this,
one can check that the values of $\hat \mu_d$
are small in comparison to the values $\hat \mu_i$ with $1 \le i<d$
and similarly that $\hat D_{id}$ and $\hat D_{dd}$
are small in comparison to
$\hat D_{ij}$ with $1 \le i,j<d$.

To make this comparison precise,
we wish to calculate scalar quantities corresponding to the
amounts of drift and diffusion both tangent and normal to the fibre.
For the one-dimensional normal direction,
the component $\hat \mu_d$ corresponds
to the normal drift and $\hat D_{dd}$ the normal component of the diffusion.
Therefore, define
$\mun(z) = |\hat \mu_d|$ and $\DDn(z) = |D_{dd}|$.

For the tangent direction, define
$\mut(z) = \sqrt{\hat \mu_1^2 + \cdots + \hat \mu_{d-1}^2}$
which is the norm of the vector
$(\hat \mu_1, \ldots, \hat \mu_{d-1})$.
Define
$\DDt(z)$ as the largest eigenvalue of the $(d-1) \times (d-1)$ submatrix
of $\hat D$.
This is equivalent to
\[
    \DDt(z) =
    \sup \{ u^\top \hat D u : \ u \in \bbR^{d-1} \times \{0\},\ \|u\|=1 \}
\]
and corresponds to the largest amount of diffusion in any direction tangent
to $\bbR^{d-1} \times \{0\}$.
Since $\hat D$ is symmetric, the value $\DDt(z)$ is also the largest singular
value of the submatrix and can easily be computed numerically
(say by the {\tt norm} function in MATLAB).
In the special case where $d = 2$, the definitions reduce down to
$\mut(z) = |\hat \mu_1|$ and $\DDt(z) = |\hat D_{11}|.$


In general, the definitions of $\mut$ and $\DDt$ depend on the choice of $A$.
If we impose the restriction that $A$ is an isometry
then one can show that $\mut$ and $\DDt$
depend only on the subspace $A \inv(\bbR^{d-1} \times 0)$
and not on the choice of $A$ itself.
Because of this property,
we only consider the case
where $A$ is an isometry.
This property holds exactly when the matrix representing $A$ is orthogonal.
The numerical methods given above for constructing this matrix ensure
that it is orthogonal.


We now summarize the algorithm.

\subsection{Local test for multiscale dynamics}

\quad 

\medskip

\textbf{\noindent Algorithm 1}

\begin{enumerate}
    \item \label{item-fullFP}
    Given an \Ito\ drift-diffusion process on a subset $\cZ$ of $\bbR^d$,
    compute the invariant density $\rho_0$ by numerically
    solving $\cL \rho_0 = 0$.

    \item \label{item-fullFPa}
    Compute the leading eigenvalues $\lam_k$ and eigenfunctions $\psi_k$
    of the adjoint operator $\cLa$.

    \item \label{item-pickfiber}
    For the eigenvalue $\lam_1  \ne  0$ and associated eigenfunction
    $\psi_1$, compute a fibre $F$ 
    which is the connected component of a level set
    $\psi_1 \inv(\{s\})$ for some $s \in \bbR$.

    \item \label{item-presample}
    Define a finite subset $Q \subset F$ consisting of points uniformly sampled
    over $F$.
    Define weights $w_z := c\, \rho_0(z)$ for $z \in Q$
    with $c > 0$ chosen such that $\sum_{z \in Q} w_z = 1$.

    \item \label{item-sample}
    For each point $z \in Q$:
    \begin{enumerate}
        \item Using linear regression on points of $F$ near $z$,
        find a hyperplane $H$ such that 
        $H+z$ is approximately tangent to $F$ at $z$.

        \item Construct an isometry $A:\bbR^d \to \bbR^d$
        taking $H$ to $\bbR^{d-1} \times \{0\}$.\\
        For a neighbourhood $U$ of $z$,
        the image $A(U \cap F)$
        is equal to a graph $\{(v,g(v)) : v \in V\}$ of a function
        $g:V \to \bbR$ with $V \subset \bbR^{d-1}$.

        \item Approximate the first and second derivatives of $g$ using
        polynomial regression on points in
        $A(F)$ near $A(z)$.

        \item Compute the drift vector $\hat \mu$ and diffusion matrix $\hat D$
        at the point $(v,g(v)) = A(z)$
        using equations \eqref{eqn-hat-mu}--\eqref{eqn-last-D}.
        From this, compute the normal and tangent components
        $\mun(z)$, $\mut(z)$, $\DDn(z)$, and $\DDt(z)$
        as explained above.
    \end{enumerate}
    \item \label{item-build-average}
    Compute the average values
    \begin{align*}
        \mun \avg = \sum_{z \in Q} w_z\, \mun(z), \quad
        \mut \avg = \sum_{z \in Q} w_z\, \mut(z),\\
        \DDn \avg = \sum_{z \in Q} w_z\, \DDn(z), \quad
        \DDt \avg = \sum_{z \in Q} w_z\, \DDt(z).
    \end{align*}
    \item \label{item-average}
    If either of $\mut \avg$ or $\DDt \avg$ is significantly larger than
    both $\mun \avg$ and $\DDn \avg$,
    this is evidence that the full system has multiscale behaviour.
\end{enumerate}
In the above algorithm, steps 3 through 7
may be repeated for a number of distinct fibres
in order to test for multiscale behaviour throughout the domain.

\section{Estimating the time-scale separation} \label{sec-timescale} 

We now propose a further test for multiscale behaviour
which also estimates the \emph{magnitude} of the time-scale
separation between the slow and fast dynamics.
This estimation is achieved by comparing the spectrum of the operator $\cL$
on the full system to the spectrum of an operator defined
by dynamics along the fibre $F$.

Recall from Section \ref{sec-generator} that the leading eigenvalues
of $\cL$ (which are the same as those of its adjoint $\cLa$)
correspond to the slowest
rates of decay for observables under the dynamics of the system.
Therefore, they give an estimate of the speed of the slow dynamics
which can be computed from the full system.

We now define a Fokker--Planck operator $\cLtan$ for the dynamics
tangent to the fast fibre $F$.
For a multiscale system, the dynamics associated to $\cLtan$ is
entirely fast and does not capture any of the slow dynamics.
Therefore, its leading eigenvalues $\lamtan_k$
correspond to the speed of the fast dynamics alone.

We are assuming that the slow dynamics evolves on a time-scale of $O(1)$
and the fast dynamics on a time-scale of $O(1/\ep)$ for some unknown
value $\ep$.
The leading eigenfunctions $\psi_k$ of $\cLa$
are observables which decay to zero at the rate of the slow dynamics
of the systems.
Because of this, the real part of the associated eigenvalues $\lam_k$ are
$O(1)$.
The leading eigenfunctions of $\cLtan^*$ are observables along the fibre $F$
which decay at a speed associated to the fast dynamics and so the
associated eigenvalues $\lamtan_k$ are $O(1/\ep)$.
The ratio of the real parts of the leading $\lam_k$
and the leading $\lamtan_k$ then gives a quantitative estimate of the time-scale
separation between the slow and fast variables.
See also the discussion in \cite{Crommelin11,FGH1}.

We now give an algorithm for determining these values $\lamtan_k$ numerically.
Let $F$ be a fibre as computed in Algorithm 1.
Further suppose that $F$ is such that it can be represented \emph{globally}
as the graph of a function.  That is,
after a linear change of coordinates $A:\bbR^d \to \bbR^d$,
the image $A(F)$ is equal to the set
$\{ (v,g(v)) : v \in V \}$ for some function $g : V \to \bbR$ defined on a subset $V$
of $\bbR^{d-1}$.
One may then define a drift-diffusion process
on $V$ where
the drift and diffusion components are given by $\hat \mu_i$
and $\hat D_{ij}$ for $1  \le  i,j  \le  d-1$ as defined in the previous section.
Using this data, the leading eigenvalues $\lamtan_k$ for the operator $\cLtan$
may be computed numerically.

We summarise the algorithm for the estimation of the time-scale separation as

\medskip

{\noindent \bf Algorithm 2A}
\begin{enumerate}
    \item \label{item-againFPa}
    Given an \Ito\ drift-diffusion process on a subset $\cZ$ of $\bbR^d$,
    compute the leading eigenvalues $\lam_k$ and eigenfunctions $\psi_k$
    of the adjoint $\cLa$ of the Fokker--Planck operator.

    \item For the eigenvalue $\lam_1  \ne  0$ and associated eigenfunction
    $\psi_1$, compute a fibre $F$ 
    which is the connected component of a level set
    $\psi_1 \inv(\{s\})$ for some $s \in \bbR$.

    \item Find a linear change of coordinates $A:\bbR^d \to \bbR^d$ such that
    $A(F)$ is equal to a graph of a function
    $g:V \to \bbR$ with $V \subset \bbR^{d-1}$.

    \item Define a drift-diffusion process on $V$ where for each point $v \in V$
    the drift is given by the first $d-1$ components of
    $\hat \mu(z)$ and
    the diffusion by the $(d-1)\times(d-1)$ submatrix of
    $\hat D(z)$
    where $z = A \inv(v,g(v))$.

    \item \label{item-VFP}
    Numerically compute the leading eigenvalues $\lamtan_k$
    of the Fokker--Planck
    operator $\cLtan$ (and/or its adjoint) for this drift-diffusion
    process on $V$.

    \item Compare the eigenvalues $\lamtan_k$ computed in step 5
    to the eigenvalues $\lam_k$ computed in step 1.
    If the ratios $\Real(\lam_k) / \Real(\lamtan_k)$ are all of similar
    magnitudes and are close to zero, then
    this implies multiscale behaviour and the ratio gives an
    estimate of the time-scale separation of the slow and fast dynamics.
\end{enumerate}
For this algorithm,
it is not necessary to estimate the partial derivatives of $g$.
This is because these derivatives do not occur in equations
\eqref{eqn-hat-mu} and \eqref{eqn-hat-D}.

Ideally, Algorithms 1 and 2A should both be performed on a given system
to test for multiscale behaviour.
Algorithm 1 tests that on a fibre $F$ the tangent components of the
drift and diffusion dominate the normal components.
Therefore, it justifies approximating the fast dynamics by
a drift-diffusion process restricted to the fibre
and defined purely by these tangent components.
Algorithm 2A then uses the eigenvalues of the resulting operator
$\cLtan$ associated to this process
in order to estimate the speed of the dynamics restricted to the fibre.


Algorithm 2A assumes that the fibre $F$ can be expressed globally as the
graph of a function.
If the shape of $F$ is such that finding a graph is not possible, then
one could instead estimate these eigenvalues by
using methods of trajectory
integration and projection to the fibre as described in \cite{FGH1}.
This alternative method does not follow the general
trajectory-free framework
proposed in this work.
It may therefore be slower in general and subject to the concerns listed in
the introduction.
In the next section, we give a variant of Algorithm 2A which is always
applicable in the specific case where the slow and fast directions
are one-dimensional.

\section{Parameterizing by arc length} \label{sec-arc} 

In the case that the fast fibre is one-dimensional
(and the full system is therefore two-dimensional),
an alternative to expressing the fibre as a graph is to parameterize the fibre
by arc length.
This is always possible, as topologically a
one-dimensional fibre $F$ will either be a line or a
circle.

Given a representation of the curve $F$, compute a sequence of
points $\{q_n\}$ such that each point is at a uniform distance from the previous.
To analyse the drift and diffusion at $q_n$,
first approximate the tangent to the curve at $q_n$
by the vector $q_{n+1} - q_{n-1}$.
Apply a rotation $A:\bbR^2 \to \bbR^2$ such that
$A(q_{n+1} - q_{n-1})$ is roughly horizontal.
Writing $(v_n,y_n) = A(q_n)$,
apply a regression method to find a polynomial $p_g$
such that $y_m \approx p_g(v_m)$ for points
$(v_m, y_m)$ near $(v_n,y_n)$.
By adjusting the angle of the rotation $A$,
one may further assume that the
first derivative $p_g'$ satisfies $p_g'(v_n) = 0$.
Using the rotation $A$ composed with the
change of coordinates mapping $(v,y)$ to $(v,y-p_g(v))$
compute the drift and diffusion data tangent and normal to the curve.
Applying \eqref{eqn-hat-mu}--\eqref{eqn-last-D} in this specific case
and using that $p_g'(v_n) = 0$, the resulting coefficients are
\[
    \hat \mu_1 = \mu^A_1, \quad
    \hat \mu_2 = \mu^A_2 - \frac{1}{2} D^A_{11} p_g'', \quad
    \hat D_{11} = D^A_{11}, \quad
    \hat D_{12} = D^A_{12}, \quad
    \hat D_{22} = D^A_{22}.
\]
For $v > v_n$ the arc length of the graph
$\{(u,p_g(u)) : u \in [v_n,v] \}$
is given by
\[
    \ell(v) = \int_{v_n}^v \sqrt{1+ [p_g'(u)]^2} \, du.
\]
Using $p_g'(v_n) = 0$,
the first and second derivatives satisfy
\[
    \frac{d \ell}{d v}\Big|_{v=v_n} = 1 \qandq
    \frac{d^2 \ell}{d v^2}\Big|_{v=v_n} = 0.
\]
If we take the derivative of $v$ with respect to itself,
then
\[
    \frac{d v}{d v}\Big|_{v=v_n} = 1 \qandq
    \frac{d^2 v}{d v^2}\Big|_{v=v_n} = 0
\]
as well.  This implies that
two different parameterizations of the curve---one by arc length
and the other by $v$---agree up to second order
at $v_n$.
Therefore, the Fokker--Planck operators coincide for the two parameterizations
with the same drift and diffusion coefficients
at the point $(v_n,y_n)$.
Using this data, one can find for any point on $F$
the components tangent to the
curve.
One can therefore define a drift-diffusion process on the curve $F$
as parametrized by arc length
and numerically compute the leading eigenvalues of the
corresponding Fokker--Planck operator as before.

We present this as a variant of Algorithm 2A to estimate the time-scale separation:

\medskip

{\noindent \bf Algorithm 2B} 
\begin{enumerate}
    \item Given an \Ito\ drift-diffusion process on a subset $\cZ$ of $\bbR^2$,
    compute the leading eigenvalues $\lam_k$ and eigenfunctions $\psi_k$
    of the adjoint $\cLa$ of the Fokker--Planck operator.

    \item For the eigenvalue $\lam_1  \ne  0$ and associated eigenfunction
    $\psi_1$, compute a fibre $F$ 
    which is the connected component of a level set
    $\psi_1 \inv(\{s\})$ for some $s \in \bbR$.

    \item Define a uniformly spaced finite sequence of points $\{q_n\}$ along $F$.

    \item Using the techniques in this section, calculate for each
    $\{q_n\}$ the components $\hat \mu_1$ and $\hat D_1$
    which are tangent to $F$.
    These values define a drift-diffusion process on $F$.

    \item Numerically compute the leading eigenvalues $\lamtan_k$
    of the Fokker--Planck
    operator $\cLtan$ (and/or its adjoint) for this drift-diffusion
    process on $F$.

    \item Compare the eigenvalues $\lamtan_k$ computed in step 5
    to the eigenvalues $\lam_k$ computed in step 1.
    If the ratios $|\lam_k| / |\lamtan_k|$ are all of similar magnitudes
    and are close to zero, then
    this implies multiscale behaviour and the ratio gives an
    estimate of the time-scale separation of the slow and fast dynamics.
\end{enumerate}

%
%
%
%
%

\section{Determining the reduced dynamics} \label{sec-reduced} 

Once the existence of multiscale dynamics has been established,
one further step is to compute the reduced slow dynamics.
We give here a simple technique
to perform the reduction to a
lower dimensional SDE.
As with the other methods developed in this paper,
the algorithm does not rely on trajectory integration.

First, consider a multiscale system where model reduction is applicable.
That is, following the formalism given in \cite{givon2004, FGH1},
there is a projection $\cP:\cZ \to \cX$ from the domain of the full
system $\cZ$ to a smaller dimensional space $\cX$
and an SDE on $\cX$ such that
orbits of the reduced system on $\cX$
resemble
projected orbits of the full system.
Our goal is to reconstruct this reduced dynamics.

As before, let $\cL$ denote the Fokker--Planck operator
for the full SDE on $\cZ$
and let $(\lam_k,\psi_k)$ be the eigenpairs for the adjoint $\cLa$.
For ease of notation, we assume here that all eigenvalues are simple.
Similar reasoning holds in the more general case.
Let $\tilde \cL$ denote the Fokker--Planck operator
for the reduced SDE on $\cX$
and let $(\tilde \lam_k, \tilde \psi_k)$ be the eigenpairs
of $\tilde \cLa$.
Assume both sequences are indexed by decreasing real part
as in \eqref{eqn-decreal}.

Recall that
the leading eigenvalues $\lam_k$ of $\cLa$
are associated to the speed of the slow dynamics and
the eigenfunctions $\psi_k$
are approximately constant along fast fibres.
In fact, using operator approximation theory,
one may show that for $\ep$ sufficiently small,
the approximations
$\lam_k \approx \tilde \lam_k$
and $\psi_k \approx \tilde \psi_k \circ \cP$ hold \cite{Crommelin11}.

Suppose $C \subset \cZ$ is a smooth curve transverse to the fast
fibre
which intersects each fibre exactly once.
Define a projection $P_c:\cZ \to C$
by letting $P_c(z)$ be the unique point on $C$ which is on the same fast fibre
as $x$.
Then $C$ may be identified with $\cX$ and $P_c$ with $\cP$.
The above approximations then imply that $(\lam_k, \psi_k|_C)$
is an approximation of the eigenpair $(\tilde \lam_k, \tilde \psi_k)$.
Moreover, the objects $\lam_k$, $\psi_k$, and $C$ can all be computed
numerically from the operator $\cLa$ for the full system.
This means that even though we have no direct way
of representing the operator
$\tilde \cLa$ associated to the reduced dynamics on $\cX$,
we still have an indirect method of computing
its leading eigenpairs.

\begin{figure}[tp]
\begin{centering}
\includegraphics{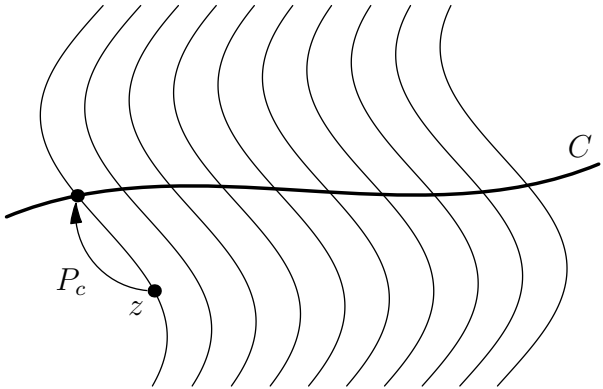}
\end{centering}
\caption{A depiction of the curve $C$, which intersects each fast fibre in a
single point and the resulting projection $P_c:\cZ \to C$
defined 
by letting $P_c(z)$ be the unique point on $C$ which is on the same fast fibre
as $z$.
We define reduced slow dynamics
on the curve $C$.
}
\label{fig-section}
\end{figure}

There are a number of techniques for reconstructing
the drift and diffusion of a system from its spectral data.
For instance, see the survey \cite{Hansen1998}
and more recent work in \cite{Crommelin06,Crommelin11}.
We use a method attributed in \cite{Hansen1998} to S.~G.~Demoura.
Let $\tilde \mu$ and $\tilde D$ be the (as yet unknown) 
coefficients of the drift and diffusion for the reduced SDE on $C$.
Consider the formula \eqref{eqn-FPa} for the adjoint of the
Fokker--Planck operator in the one-dimensional setting
along the curve $C$.
One sees that each eigenpair $(\tilde \lam_k, \tilde \psi_k)$ of $\tilde \cLa$
satisfies the equation
\begin{equation} \label{eqn-demoura}
    \tilde \lam_k  \tilde \psi_k = \tilde \mu \tilde \psi_k' + \tfrac{1}{2} \tilde D \tilde \psi_k''
\end{equation}
where the prime denotes derivatives along the curve $C$.
Using $\psi_k|_C$ as an approximation for $\tilde \psi_k$,
one may estimate the derivatives $\tilde \psi_k'$ and $\tilde \psi_k''$
numerically.
By considering \eqref{eqn-demoura} for
two distinct eigenpairs, say
$(\tilde \lam_1, \tilde \psi_1)$
and
$(\tilde \lam_2, \tilde \psi_2)$,
one may solve the linear system of equations
\begin{equation} \label{eqn-demoura-mat}
        \begin{bmatrix}
        \tilde \psi_1' & \tfrac{1}{2} \tilde \psi_1'' \\
        \tilde \psi_2' & \tfrac{1}{2} \tilde \psi_2''  \end{bmatrix}
        \begin{bmatrix}
        \tilde \mu \\
        \tilde D  \end{bmatrix}
    =
        \begin{bmatrix}
        \tilde \lam_1 \tilde \psi_1 \\
        \tilde \lam_2 \tilde \psi_2  \end{bmatrix}
      \end{equation}
pointwise along $C$ in order to solve
for the drift $\tilde \mu$ and diffusion $\tilde D$
of the reduced system.

This technique requires an accurate estimate
of the first and second derivatives of the eigenfunctions.
In some settings, accurate derivative estimation
is not possible and more sophisticated
techniques must be applied, such as 
solving a quadratic programming problem
based on the eigenfunction equations
\cite{Gobet2004, Crommelin06, Crommelin11}.

\section{Numerical example} \label{sec-numerics} 

We now apply these techniques to an example studied
in \cite{Crommelin06} and given by

\begin{align} \label{eqn-original-x}
    \dot x &= \sin y + \sqrt{1 + \frac{1}{2} \sin y}\ \dot W_x,\\
    \label{eqn-original-y}
    \dot y &= \frac{1}{\ep}\bigl[- y + \sin x \bigr] + \frac{1}{\sqrt{\ep}}\ \dot W_y
\end{align}
where $W_x$ and $W_y$ are independent Wiener processes.
This SDE is
defined on $[0,2 \pi] \times \bbR$ with periodic boundary conditions on $x$.
In order to have an example with nonlinear fast fibres,
we consider the system after a change of coordinates taking
$(x, y)$ to $(x+\sin(y), y)$.
After this change, the system is given by the
(admittedly much uglier)
equations

\begin{align} \label{the-system-start}
    \dot x =& \left[
             \sin y +
             \frac{\cos y}{\ep}\bigl(\sin(x-\sin y) - y \bigr)
             - \frac{\sin y}{2 \ep}
             \right]\\
             \nonumber
             &+
             \sqrt{1 + \frac{1}{2}\sin y}\ \dot W_x
             +
             \frac{1}{\sqrt{\ep}} \cos y\ \dot W_y, \\
    \label{the-system-end}
    \dot y =& \frac{1}{\ep}\bigl[\sin(x - \sin y) - y \bigr]
              + \frac{1}{\sqrt{\ep}}\ \dot W_y.
\end{align}
Note that terms involving $\ep$ now occur in the formulas for both $\dot x$
and $\dot y$ and so
these equations have no obvious slow-fast splitting.
As in \cite{Crommelin06} we fix $\ep = 10^{-3}$.
We now apply the steps of Algorithm 1 to test for multiscale behaviour.

For step 1 of the algorithm,
we compute the invariant density $\rho_0$
by solving numerically for the leading eigenfunction of
the Fokker--Planck operator $\cL$
using a collocation method
\cite{Boyd01}.
As $x$ is periodic and $y$ is not,
we restrict the domain to $[0, 2 \pi] \times [-L, L]$
with $L = \Lmat$
and use as basis the tensor product of
a Fourier basis in $x$ and
a Chebyshev basis in $y$.
A $\Nmat \times \Nmat$ grid of points is used for the collocation.
Since one expects the density to decay to zero as $|y| \to \infty$,
we impose a Dirichlet boundary condition that
$\rho_0$ is zero on
$[0, 2 \pi] \times \{-L\}$ and $[0, 2 \pi] \times \{L\}$.

For step 2, we solve for the leading eigenpairs
$(\lam_k, \psi_k)$ of the adjoint operator $\cLa$
using the same domain and basis.
However, since the eigenfunctions of $\cLa$ corresponding to the eigenvalue $\lam_0 = 0$ are given by the constant functions,
we impose the boundary condition
$\frac{\partial \psi}{\partial y} = 0$ on
$[0, 2 \pi] \times \{-L\}$ and $[0, 2 \pi] \times \{L\}$.
\begin{figure}[tp]
\begin{centering}
\includegraphics{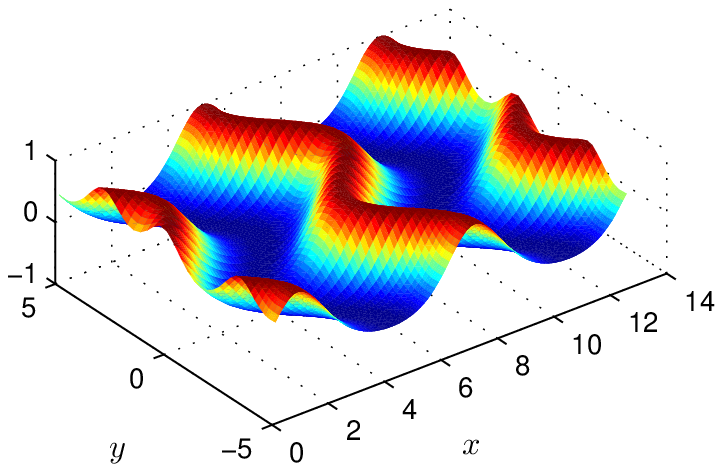}
\includegraphics{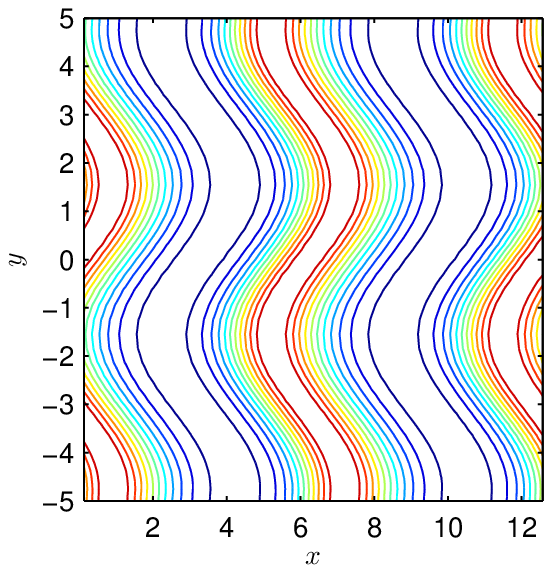}
\end{centering}
\caption{The real part of the eigenfunction $\psi_1$
corresponding to
the eigenvalue $\lam_1 = \lammat$
for the adjoint $\cLa$ of the
Fokker--Planck operator for the SDE given by
\eqref{the-system-start}--\eqref{the-system-end}
with $\ep = 10^{-3}$.
The function is plotted for $x \in [0, 4 \pi]$
using the $2 \pi$-periodicity of the $x$ variable.
Also plotted are the resulting level sets,
used as a numerical approximation of the fast fibres.
}
\label{fig-eigsurf}
\end{figure}
The leading eigenvalues are given in
the first column of Table \ref{table-eigs}.
After the eigenvalue $\lam_0 = 0$
(numerically computed as $\lamone$),
the eigenvalues with greatest real
part are 
$\lam_1 = \lammat$ and its conjugate $\lam_2 = \bar \lam_1$.
The eigenfunction $\psi_1$
associated to $\lam_1$ is complex valued
and the real part of this function
is plotted in Figure \ref{fig-eigsurf}.


This real part is used to compute level sets which
closely approximate fast fibres of the system
when multiscale behaviour is present (cf. step 3).
For concreteness, we use a fibre $F$
defined by the level set of $\Real \psi_1$
which passes through the point
$(x,y) = (5,0)$.
Other choices of level sets yield similar results.

Once this fibre $F$ is computed,
we define a finite subset $Q$ by sampling points
along $F$ where each point is at a distance
$0.1$ from the previous.
Taking weights $w_z$ proportional to $\rho_0(z)$,
this defines a measure on the fibre (cf. step 4).

\begin{figure}[tp]
\begin{centering}
\includegraphics{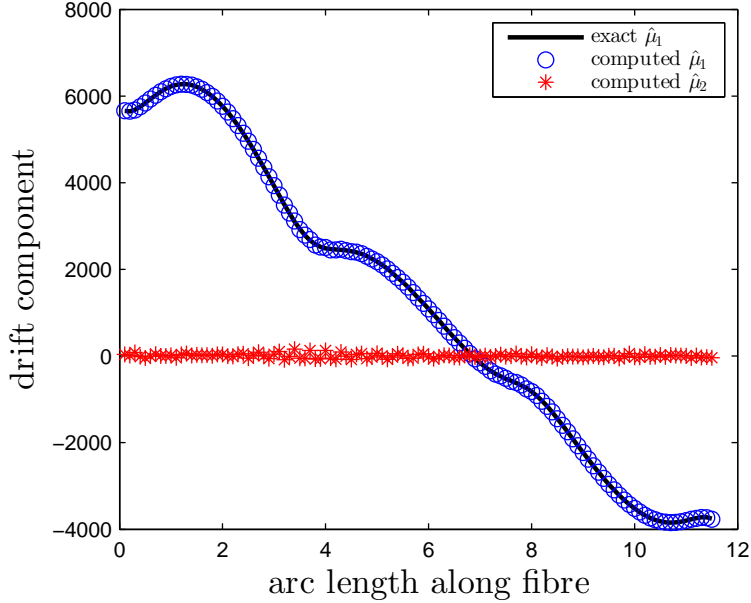}
\end{centering}
\caption{The drift components computed in coordinates adapted to the
fibre. The $\hat \mu_1$ component corresponds to drift along the fibre
and $\hat \mu_2$ to drift normal to the fibre.
}
\label{fig-fiberdrift}
\end{figure}
\begin{figure}[tp]
\begin{centering}
\includegraphics{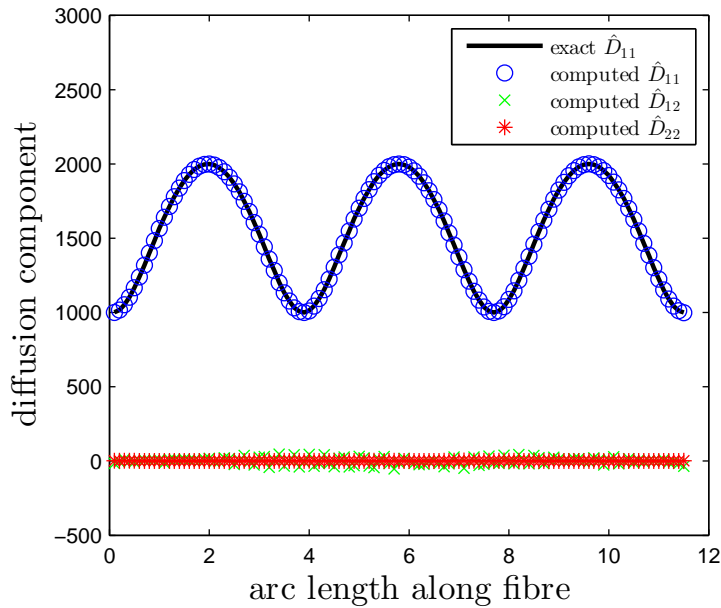}
\end{centering}
\caption{The diffusion components computed in coordinates adapted to
the fibre. The $\hat D_{11}$ component corresponds to diffusion along the fibre.
}
\label{fig-fiberdiff}
\end{figure}

For each point in $Q$, we construct a change of coordinates
(cf. step 5).
This change of coordinates yields
at each point a new drift vector $\hat \mu$
and diffusion matrix $\hat D$.
The components of the computed $\hat \mu$ and $\hat D$ are plotted
in Figures \ref{fig-fiberdrift} and \ref{fig-fiberdiff}.
Note that the components $\hat \mu_1$ and $\hat D_{11}$
are the largest in magnitude,
and these correspond to the dynamics tangent to the fibre $F$.

From the original equations
\eqref{eqn-original-x} and \eqref{eqn-original-y},
one sees that the dynamics on a fast vertical fibre
in these original coordinates
are given by a process with drift
$(- y + \sin x)/\ep$
and diffusion ${1}/{\sqrt{\ep}}.$
These drift and diffusion formulas can be transformed under the change
of coordinates $(x,y) \mapsto (x+\sin(y),y)$
to give exact analytical formulas for the
drift $\hat \mu_1$ and diffusion $\hat D_{11}$ along the fibre
when it is parameterized by arc length.
These exact functions for $\hat \mu_1$ and $\hat D_{11}$ are also plotted in
Figures \ref{fig-fiberdrift} and \ref{fig-fiberdiff}
and agree closely with the computed values.
Since these analytical formulas are long and not elucidating,
we do not
include them here.




To compare the tangent and normal components
quantitatively,
we apply step 6 of Algorithm 1,
computing the weighted averages over the fibre.
The resulting averages for the drift are
$\mut \avg = \avgdatafour$ and
$\mun \avg = \avgdatafive$
and for the diffusion are
$\DDt \avg = \avgdataone$
and $\DDn = \avgdatathree$.
Note that $\mut \avg$ and $\DDt \avg$ are both significant larger
than $\mun \avg$ and $\DDn \avg$.
By step 7, we conclude that the system has multiscale behaviour.

%

\bigskip

We now apply Algorithms 2A and 2B to the system.
Before discussing the numerics, first note that
from the original equations
\eqref{eqn-original-x} and \eqref{eqn-original-y},
one sees that the dynamics on a fast vertical fibre
in these original coordinates
is given by an Ornstein--Uhlenbeck process with drift
$(-y+\sin x)/\ep$
and diffusion ${1}/{\sqrt{\ep}}.$
For such processes, the eigenvalues of the Fokker--Planck operator
can been determined analytically
\cite{sjogren-survey}.
For \eqref{eqn-original-y}, the eigenvalues are
given exactly by $\lambda_k=-k/\ep$
for integers $k \ge 0$.
Since the spectrum of the operator is not affected by a change of coordinates,
these are also the exact values for the eigenvalues of the operator $\cLtan$
associated to the fast dynamics of the system given by
\eqref{the-system-start}--\eqref{the-system-end}.

In the application of Algorithm 1 above,
we computed the values
$\hat \mu_1$ and $\hat D_{11}$ at points uniformly distributed along the fast
fibre $F$.
These define a drift diffusion process along $F$
with associated Fokker--Planck operator $\cLtan$.
Applying Algorithm 2B,
we calculate leading eigenvalues $\lamtan_k$ of this operator
$\cLtan$ numerically.
A fourth-order central finite difference method was used.
The computed eigenvalues are given in the second column of Table
\ref{table-eigs}.
Note that the eigenvalues $\lamtan_k$ are all of order
$10^3$ and the eigenvalues $\lam_k$
are of order $1$.
Applying the final step of Algorithm 2B,
the ratios
$\Real(\lam_k)/\Real(\lamtan_k)$ are all between $10^{-4}$
and $10^{-3}$
and give an estimate of the time-scale separation.
This agrees with the value
$\ep = 10^{-3}$ used in defining the system.

For each computed eigenvalue $\lamtan_k$, the relative error
\begin{equation} \label{eq-err}
    \err_k = \frac{|\lamtan_k - (\ideallam)|}{|\ideallam|}
      \end{equation}
is also given in the table.

\medskip

We now consider Algorithm 2A
and perform eigenvalue calculations by
expressing the fibre globally as a graph.
A fast fibre $F$ is of the form $\{(x,y) : x = \sin(y) + c\}$
for some constant $c$.  Therefore, if $A$ is a rotation matrix
\[
    A =
        \begin{bmatrix}
        \cos(\theta) & -\sin(\theta) \\
        \sin(\theta) & \cos(\theta)  \end{bmatrix}
\]
with $\theta$ between $45$ degrees and $135$ degrees,
then $A(F)$ is the graph of a function $g:\bbR \to \bbR$.
For values of $\theta$ at every $5$ degrees between $\minangle$ and $\maxangle$,
we applied Algorithm 2A with the corresponding $A$
and
computed seven leading eigenvalues
$\lamtan^{(\theta)}_0 > \lamtan^{(\theta)}_1 > \cdots > \lamtan^{(\theta)}_6$
for the associated Fokker--Planck operator.
This computation was achieved by first determining the values $\hat\mu_1(z_n)$
and $\hat D_{11}(z_n)$ for points $z_n := (v_n, g(v_n)) \in A(F)$ 
where the $v_n$ are given by $200$ evenly spaced numbers
across an interval in $\bbR$.
Then, these sampled values of $\hat\mu_1$ and $\hat D_{11}$
are used to approximate the Fokker--Planck operator
$\cLtan$ and the eigenvalues are solved using a fourth order
central finite difference method.
As noted above,
the leading eigenvalues for the actual system are given by $\ideallam$
and
the computed values $\lamtan^{(\theta)}_k$ agreed with these theoretical values
with a relative error
less than one percent for all computed $\theta$
and $1  \le  k  \le  6$.
This shows that Algorithm 2A produces estimates of the eigenvalues
$\lamtan_k^{(\theta)}$ for the fibre dynamics that are robust
with respect to the parameterization the fibre.

\medskip

Algorithms 1 and 2
have established the presence of multiscale behaviour for the system
and determined the time-scale separation.
As a final step,
we compute the drift and diffusion data for the reduced slow system
using the technique described in Section \ref{sec-reduced}.
We restrict the leading eigenfunctions $\psi_1$ and $\psi_2$
of the full system to the line $C = [0,2 \pi] \times \{0\}$ which is transverse to
the fast fibres.
Each eigenfunction $\psi_k$ on the full space was solved numerically
by collocation
and the result is represented by the values that $\psi_k$ takes on a
$\Nmat \times \Nmat$ grid on $[0, 2 \pi] \times [-L, L]$.
Restriction of this eigenfunction to $C$ may therefore be performed
simply by taking the appropriate row from the matrix representing
the computed function.
The restriction $\psi_k|_C$ is then given
by the values that the function takes on a set $X \times \{0\}$
consisting of $\Nmat$ points uniformly spaced along $C$.
At each point $x \in X$, we look at a window of points consisting of
$x$ and the ten other points in $X$ closest to
$x$. Using these eleven data points, we fit a cubic polynomial to $\psi_k|_C$.
This polynomial is used to estimate
the first and second derivatives of $\tilde \psi_k \approx \psi_k|_C$ at
the point $x$
which
are then used to 
solve pointwise the linear equations given by \eqref{eqn-demoura-mat}.
\begin{table}[tp]
\caption{Leading eigenvalues $\lam_k$ computed
for the adjoint of the
Fokker--Planck operator for the SDE given by
\eqref{the-system-start}--\eqref{the-system-end}
with $\ep = 10^{-3}$.
Also computed are the leading eigenvalues $\lamtan_k$ of
the Fokker--Planck operator $\cLtan$ on a fast fibre $F$.
The ratios 
{$\Real(\lam_k)/\Real(\lamtan_k)$}
give a computable estimate
of the time-scale separation $\ep$.
Analytically, $\lamtan_k$ should be equal to $\ideallam = -k \times 10^3$
and the relative error in computation, as defined in \eqref{eq-err}, is also given.
}
\centering
\begin{tabular}{c c c c c}
\toprule
\vspace{1mm}
{$k$} &
{$\lam_k$} &
{$\lamtan_k$} &
{$\Real(\lam_k)/\Real(\lamtan_k)$} &
{$\err_k$}
\\
\hline
\ensuremath{0} & \ensuremath{2.0971 \times 10^{-11}} & \ensuremath{-3.3707 \times 10^{-6}} & --- & --- \\
\ensuremath{1} & \ensuremath{-0.6467 + 0.1097 i} & \ensuremath{-9.9275 \times 10^{2}} & \ensuremath{6.5139 \times 10^{-4}} & \ensuremath{7.2535 \times 10^{-3}} \\
\ensuremath{2} & \ensuremath{-0.6467 - 0.1097 i} & \ensuremath{-2.0325 \times 10^{3}} & \ensuremath{3.1816 \times 10^{-4}} & \ensuremath{1.6268 \times 10^{-2}} \\
\ensuremath{3} & \ensuremath{-2.0508 + 0.2465 i} & \ensuremath{-3.0931 \times 10^{3}} & \ensuremath{6.6303 \times 10^{-4}} & \ensuremath{3.1044 \times 10^{-2}} \\
\ensuremath{4} & \ensuremath{-2.0508 - 0.2465 i} & \ensuremath{-4.0126 \times 10^{3}} & \ensuremath{5.1110 \times 10^{-4}} & \ensuremath{3.1492 \times 10^{-3}} \\
\ensuremath{5} & \ensuremath{-4.4543 + 0.3912 i} & \ensuremath{-4.9614 \times 10^{3}} & \ensuremath{8.9780 \times 10^{-4}} & \ensuremath{7.7173 \times 10^{-3}} \\
\ensuremath{6} & \ensuremath{-4.4543 - 0.3912 i} & \ensuremath{-5.9946 \times 10^{3}} & \ensuremath{7.4306 \times 10^{-4}} & \ensuremath{8.9551 \times 10^{-4}} \\

\bottomrule
\end{tabular}
\label{table-eigs}
\end{table}
\begin{figure}[tp]
\begin{centering}
\includegraphics{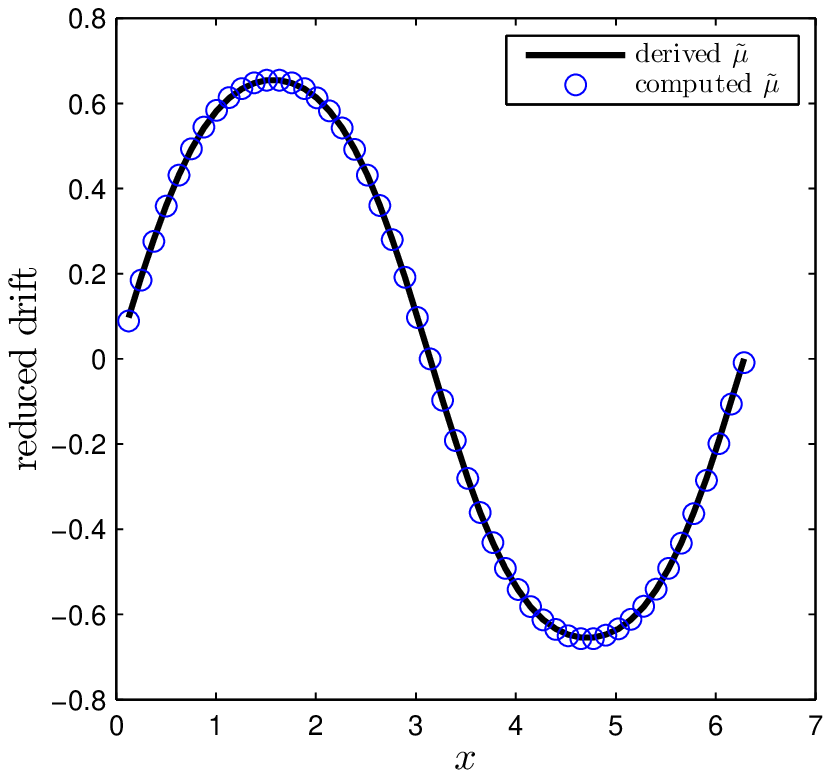}
\end{centering}
\caption{The drift component $\tilde \mu$ computed for the reduced slow dynamics
along the line $C = [0, 2 \pi] \times \{0\}$.
Plotted as circles is the value determined by computation.
Plotted as a solid line is the value derived by
homogenization theory in the limit $\ep\to 0$.
}
\label{fig-slowdrift}
\end{figure}
\begin{figure}[tp]
\begin{centering}
\includegraphics{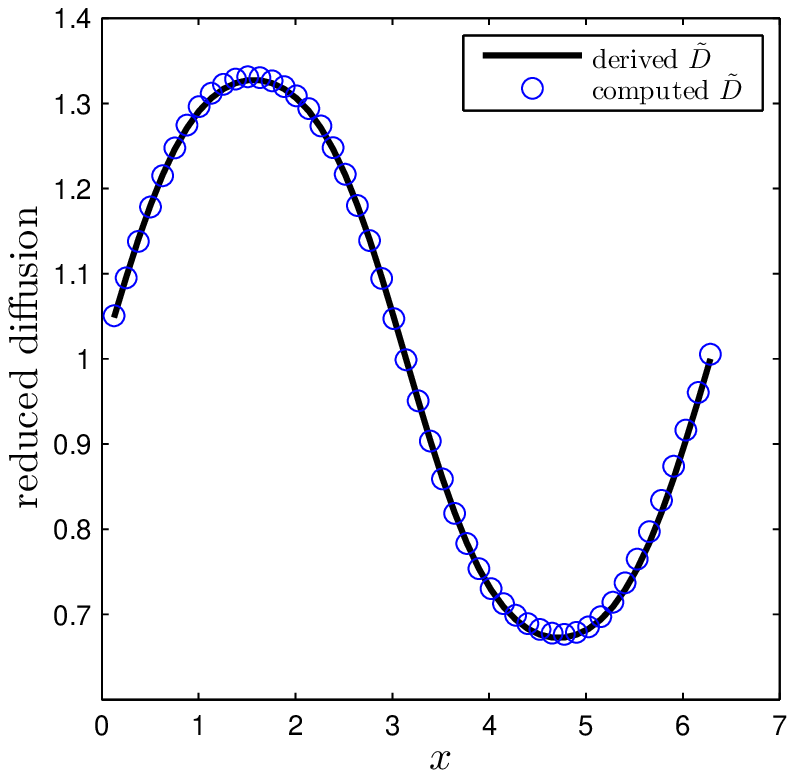}
\end{centering}
\caption{The diffusion component $\tilde D$ computed for the reduced slow
dynamics along the line $C = [0, 2 \pi] \times \{0\}$.
Plotted as circles is the value determined by computation.
Plotted as a solid line is the value derived by
homogenization theory in the limit $\ep\to 0$.
}
\label{fig-slowdiff}
\end{figure}
The computed values
for the drift $\tilde \mu$ and diffusion $\tilde D$
of the reduced system
are plotted in Figures \ref{fig-slowdrift} and
\ref{fig-slowdiff}.
Using homogenization techniques (see \cite{Crommelin06}) the reduced slow system of the multi-scale system \eqref{eqn-original-x}--\eqref{eqn-original-y} can be determined as
\begin{align}
dX = e^{-\frac{1}{4}} \sin(\sin(X)) \, dt + \left(1 + \tfrac{1}{2}e^{-\frac{1}{4}} \sin(\sin(X))\right) dW_t\; .
\end{align}
The drift and diffusion functions of this reduced slow equation are also plotted in the figures and the values determined by
computation agree closely with those given by homogenization theory.

\section{Discussion} \label{sec-discussion} 

To conclude, we look at the computational overhead of the techniques introduced
in this paper in comparison to other methods for analyzing multiscale systems.
In the example in Section \ref{sec-numerics},
building the matrices used
in steps 1 and 2 of Algorithm 1
required evaluating the functions for the drift and diffusion 
at every point of a
$\Nmat \times \Nmat$ grid.
Other steps in Algorithms 1, 2A, and 2B
consider the drift and diffusion
on lower-dimensional
fibres and require even fewer evaluations.
Overall, the computation of the drift and diffusion data at various points
do not significantly contribute to the run time of the algorithms.
In fact, the construction of the matrices representing the partial
derivatives are more time consuming.
Since there were no severe bottlenecks,
all of the computations in Section \ref{sec-numerics} were performed with a
total runtime of
25 seconds on a desktop computer.\footnote{
The computer used was an Intel Core i7-4770 with four CPUs at 3.40 GHz
and 16 GB of RAM running Ubuntu 12.04 with 64-bit Linux kernel
3.13.0-39-generic and MATLAB R2013b.}
This running time is small in comparison
to techniques which rely on trajectory integration to analyse
a system.

A time series long enough to accurately capture the statistics of the slow dynamics
of 
\eqref{eqn-original-x}--\eqref{eqn-original-y}
could easily be of length
$10^6$ or longer.
Assuming the time step used to compute this orbit at least on the order of
the time-scale separation $\ep \inv = 10^3$, this corresponds to
$10^{9}$ individual time steps of the Euler--Maruyama method.

This paper is a continuation of the work developed in \cite{FGH1}.
There, we relied solely on trajectory integration to analyze a system and used
Ulam's method
to compute the eigenfunctions
of $e^{t \cL}$ for a dynamical system.
To apply Ulam's method to the multiscale system,
the length of orbits computed
need to be at least on the order of the slow
dynamics and many individual segments of orbit need to be computed.
For instance, for the numerical examples given in \cite{FGH1},
each orbit segment was computed using $2 \times 10^4$ time steps
and $10^4$ individual orbits were computed in
each square of a $200 \times 200$ grid.
This gives a total count of $8 \times 10^{12}$ steps of the Euler-Maruyama
method where each step must compute the values of the drift and diffusion
coefficients of at point in the phase space.
In \cite{FGH1}, this computation was achieved
by splitting the construction of the matrix over many computers in
parallel.
Methods based on the infinitesimal generator avoid these issues
and can thus lead to large gains in
computation speed.

That said, the inifinitesimal generator method is not always applicable.  If
the diffusion in the actual reduced slow process is very small,
it may be dominated by
numerical diffusion in the full space arising when
computing the eigenfunctions of the
Fokker-Planck operator.
In this case, other less efficient methods,
such as Ulam's method with a large
enough time step, must be used to yield accurate results.


\bibliographystyle{plain}
\bibliography{multi}

\end{document}